\documentclass[12pt]{amsart}
\usepackage{amssymb}
\usepackage{amsfonts}
\usepackage{amsmath}
\usepackage{graphicx}

\renewcommand\appendix{\par
\setcounter{section}{0}\setcounter{subsection}{0}\setcounter{table}{0}
\setcounter{figure}{0}
\gdef\thetable{\Alph{table}}
\gdef\thefigure{\Alph{figure}}
\section*{Appendix. Quasianalytic Shifts}
\gdef\thesection{\Alph{section}}
\setcounter{section}{1}}
\newtheorem{theorem}{Theorem}
\theoremstyle{plain}
\newtheorem{lemma}[theorem]{Lemma}
\newtheorem{proposition}[theorem]{Proposition}
\newtheorem{corollary}[theorem]{Corollary}

\newtheorem{remark}[theorem]{Remark}

\numberwithin{equation}{section}

\input{tcilatex}

\begin{document}
\title[Stable Invariant Subspaces]{Stable and Norm-stable Invariant Subspaces%
}
\author{Alexander Borichev}
\address{Centre de Math\'ematiques et Informatique, Universit\'e
d'Aix-Mar\-seille I, 39 rue Fr\'ed\'eric Joliot-Curie, 13453 Marseille,
France}
\email{borichev@cmi.univ-mrs.fr}
\author{Don Hadwin}
\address{Mathematics Department\\
University of New Hampshire\\
Durham, NH 03824\\
USA }
\email{don@unh.edu}
\author{Hassan Yousefi}
\address{Mathematics Department\\
California State University\\
Fullerton, CA 92831\\
USA }
\email{hyousefi@fullerton.edu}
\thanks{The first author was partially supported by the ANR projects DYNOP
and FRAB}

\begin{abstract}
We prove that if $T$ is an operator on an infinite-dimensional Hilbert space
whose spectrum and essential spectrum are both connected and whose \textit{%
Fredholm index} is only $0$ or $1$, then the only nontrivial \textit{%
norm-stable invariant subspaces} of $T$ are the finite-dimensional ones. We
also characterize norm-stable invariant subspaces of any weighted unilateral
shift operator. We show that \textit{quasianalytic shift operators} are
points of norm continuity of the lattice of the invariant subspaces. We also
provide a necessary condition for strongly stable invariant subspaces for
certain operators.
\end{abstract}

\maketitle

%\urladdr{http://www.math.unh.edu/\symbol{126}don}

%This paper is in final form and no version of it will be submitted for publication elsewhere.

\renewcommand{\subjclassname}{\textup{2000} Mathematics Subject
Classification}

\section{Introduction}

In this paper we continue work started in \cite{CH} on stable invariant
subspaces of Hilbert-space operators. In Section 2, we show that if $T$ is
an operator on a separable Hilbert space whose Fredholm index at every point
in its semi-Fredholm domain is either $0$ or $1$ and whose spectrum and
essential spectrum are connected, then the proper norm stable invariant
subspaces of $T$ must be finite-dimensional (Theorem \ref{c13}). As a
consequence, in Section 3, we completely characterize the norm-stable
invariant subspaces of every weighted unilateral shift operator (Theorem \ref%
{wtshift}). In Section 4 we prove a semicontinuity result (Lemma \ref%
{semicontinuous}) for an index invariant (see, for instance, \cite{R}), and
we use it to provide a necessary condition for invariant subspaces of
certain operators to be (strongly) stable. Several open questions are given
in Section 5.

In our work, we need a variant of Y.~Domar's result answering a problem of
A. Shields (\cite[Problem 17]{Shi}, 1974, see also the updated edition of 
\cite{Shi} published in 1979): the nonzero invariant subspaces of a
quasianalytic weighted unilateral shift operator all have finite
co-dimension (Theorem \ref{Shields} in Appendix).

Throughout this paper, ${\mathcal{H}}$ is a separable infinite-dimensional
Hilbert spa\-ce, ${\mathcal{B}}({\mathcal{H}})$ is the set of all (bounded,
linear) operators on ${\mathcal{H}}$, and ${\mathcal{K}}({\mathcal{H}})$ is
the space of compact operators on ${\mathcal{H}}$. A closed linear subspace $%
M$ of ${\mathcal{H}}$ will be identified with the orthogonal projection $P_M$
onto $M$. If $T\in{\mathcal{B}}({\mathcal{H}})$, then $\func{Lat\,}(T)$
denotes the set of all closed invariant subspaces of $T$; alternatively, $%
\func{Lat\,}(T)$ is the set of all projections $P$ in ${\mathcal{B}}({%
\mathcal{H}})$ such that $(1-P)TP=0$. We denote the spectrum of $T$ by $%
\sigma(T)$ and the spectral radius of $T$ by $r(T)$. The image of an
operator $T\in{\mathcal{B}}({\mathcal{H}})$ in the Calkin algebra ${\mathcal{%
B}}({\mathcal{H}})/{\mathcal{K}}({\mathcal{H}})$ is denoted by $\widetilde{T}
$ and $\sigma_e(T)$ denotes the essential spectrum of $T$, which by
definition, equals $\sigma(\widetilde{T})$. Recall that $T\in{\mathcal{B}}({%
\mathcal{H}})$ is a \emph{semi-Fredholm operator }if $\mathrm{{Range}\,}(T)$
is closed and either $\dim(\mathrm{{Ker}\,}(T))$ or $\dim(\mathrm{{Ker}\,}%
(T^{\ast}))$ is finite. In this case, the \emph{index} of $T$ is defined by 
\begin{equation*}
\func{ind\,}(T)=\dim(\mathrm{{Ker}\,}(T))-\dim(\mathrm{{Ker}\,}(T^{\ast})). 
\end{equation*}
The \emph{semi-Fredholm domain} of $T$ is denoted by $\rho_{S\text{-}F}(T)$
and is defined by 
\begin{equation*}
\rho_{S\text{-}F}(T)=\{\lambda\in\mathbb{C}:\lambda-T\text{ is semi-Fredholm}%
\}. 
\end{equation*}
If $\func{ind\,}(\lambda-T)\neq\pm\infty$ for every $\lambda\in\rho_{S\text{-%
}F}(T)$, then 
\begin{equation*}
\mathbb{C}\backslash\rho_{S\text{-}F}(T)=\sigma_{e}(T). 
\end{equation*}
Recall that an operator $T\in{\mathcal{B}}({\mathcal{H}})$ is \emph{%
biquasitriangular} if and only if $\func{ind\,}(T-\lambda)=0$ for every $%
\lambda\in\rho_{S\text{-}F}(T)$. The \emph{similarity orbit} of $T$ is
defined to be 
\begin{equation*}
{\mathcal{S}}(T)=\{ATA^{-1}:A\in{\mathcal{B}}({\mathcal{H}})\text{ is
invertible}\} . 
\end{equation*}

One of the deepest results in approximation theory for operators is the 
\emph{Similarity Theorem} \cite[Theorem 9.2]{AFHV}. For our purposes we need
only the special case taken from \cite{AHV}, which we state here.

\begin{proposition}
\label{sim}\cite{AHV} Suppose $B,A\in {\mathcal{B}}({\mathcal{H}})$.

\begin{enumerate}
\item If $A$ is a normal operator and $\sigma(A)\subset\mathbb{C}%
\backslash\rho_{S\text{-}F}(B)$, then ${\mathcal{S}}(B)$ and ${\mathcal{S}}%
(A\oplus B)$ have the same norm closures.

\item If

\begin{enumerate}
\item $\sigma_{e}(B)$ has no isolated points,

\item either $\sigma_{p}(B)=\emptyset$ or $\sigma_p(B^{\ast})=\emptyset$,

\item each component of $\sigma_e (A)$ meets $\sigma_e (B)$,

\item $\sigma_{e}(B)\subset\sigma_e (A)$,

\item $\rho_{S\text{-}F}(A)\subset\rho_{S\text{-}F}(B)$ and $\func{ind\,}%
(A-\lambda)=\func{ind\,}(B-\lambda)\neq\pm\infty$ for every $%
\lambda\in\rho_{S\text{-}F}(A)$.
\end{enumerate}

Then $A$ is in the norm closure of ${\mathcal{S}}(B)$.
\end{enumerate}
\end{proposition}

An invariant subspace $P\in\func{Lat\,}(T)$ is called (strongly) \emph{stable%
} if, whenever there is a sequence $\{T_n\}$ in ${\mathcal{B}}({\mathcal{H}})
$ such that $\Vert T_n-T\Vert \rightarrow0$, there is a sequence $\{P_n\}$
with $P_n\in\func{Lat\,}(T_n)$ for $n\ge1$ such that $P_n\rightarrow P$ in
the strong operator topology (SOT). We say that $P$ is \emph{norm stable} if
we can always choose $\{P_n\}$ so that $\Vert P_n-P\Vert \rightarrow0$. We
let $\func{Lat\,}_{s}(T)$ be the collection of stable invariant subspaces of 
$T$ and $\func{Lat\,}_{ns}(T)$ be the collection of norm-stable invariant
subspaces of $T$. It is clear that $\func{Lat\,}_{ns}(T)\subset\func{Lat\,}%
_{s}(T)$ and it is easy to show that $\func{Lat\,}_{ns}(T)$ contains $\{0\}$
and ${\mathcal{H}}$. It is also easy to show that $\func{Lat\,}_{ns}(T)$ is
norm closed and that $\func{Lat\,}_{s}(T)$ is SOT-closed (see \cite{CH}).
The following question was posed in \cite{CH}:

\bigskip

\textbf{Question}: Is $\func{Lat\,}_{s}(T)$ always the SOT-closure of $\func{%
Lat\,}_{ns}(T)$?

\bigskip

J. Conway and D. Hadwin \cite{CH} gave an affirmative answer to this
question when $T$ is normal or an unweighted unilateral shift of finite
multiplicity.

In the finite-dimensional setting, the stable invariant subspaces of an
operator were characterized in \cite{BGK}, \cite{CD}, and \cite{AFHV}. In 
\cite{AFS} C. Apostol, C. Foia\c{s} and N. Salinas showed that if $T$ is a
normal operator, then the projections in $\func{Lat\,}_{ns}(T)$ are
precisely the spectral subspaces corresponding to clopen subsets of $%
\sigma(T)$. Moreover, they proved that a quasitriangular operator with
connected spectrum has no nontrivial norm-stable invariant subspace (i.e., $%
\func{Lat\,}_{ns}(T)=\{0,1\}$). The question for $\func{Lat\,}_{s}(T)$ is
much more delicate and is related to the invariant subspace problem \cite{CH}%
.

\section{Results on Lat$_{ns}.$}

We begin with a topological lemma for compact subsets of the plane. If $%
K\subset\mathbb{C}$ and $\varepsilon>0$, then we define $K_{\varepsilon}=\{z%
\in\mathbb{C}:\mathrm{{dist}\,}(z,K) <\varepsilon\}$.

We use the notation ${\mathbb{D}}=\{z\in\mathbb{C}:|z|<1\}$, $r{\mathbb{D}}%
=\{z\in\mathbb{C}:|z|<r\}$, ${\mathbb{T}}=\partial{\mathbb{D}}$.

\begin{lemma}
\label{c9} Suppose that $K$ is a nonempty compact connected subset of $%
\mathbb{C}$, $\delta>0$, and $\mathcal{U}$ is a collection of bounded
connected components of $\mathbb{C}\backslash K$ such that each $U\in%
\mathcal{U}$ contains a point whose distance to $K$ is greater than $\delta$%
. Let $V={\displaystyle\bigcup_{U\in\mathcal{U}}}U$. Then there is an $r>1$
and a univalent \textrm{(}analytic\textrm{)} function $\varphi$ on $r{%
\mathbb{D}}$ such that

\begin{enumerate}
\item $\varphi({\mathbb{T}})\subset K_{\delta}$, and

\item $V\backslash K_{\delta}\subset\varphi({\mathbb{D}})\subset
K_{\delta}\cup V$.
\end{enumerate}
\end{lemma}

\begin{proof}
If $\mathcal{U}=\varnothing$, then $V=\varnothing$, and we can choose an $%
a\in K$ and define $\varphi(z)=a+\frac{\delta}{r}z$. Now suppose that $%
\mathcal{U}\neq\varnothing$. Let $V_{1},\ldots,V_{n}$ be the elements of $%
\mathcal{U}$ that contain a closed disc of radius $\delta$. Since $K$ is
connected, each $V_{j}$ is simply connected, since it has a connected
complement. Thus there is a univalent function $f_{j}$ from the unit disc ${%
\mathbb{D}}$ to $V_{j}$ for $1\le j\le n$. Since $\{f_{j}(r{\mathbb{D}}%
):0<r<1\}$ is an open cover of $V_{j}$, there is a $t$, $0<t<1$, such that 
\begin{equation*}
\{z\in V_{j}:\mathrm{{dist}\,}(z,K)\ge\delta/2\}\subset f_{j}(t{\mathbb{D}}%
)\equiv W_{j} 
\end{equation*}
for $1\le j\le n$. Since 
\begin{equation*}
\Omega=K_{\delta}\cup\bigl({\bigcup_{1\le j\le n}}W_j\bigr)  
\end{equation*}
is an open and connected set containing the disjoint closed contractible
"discs" $\overline{W}_{1},\ldots,\overline{W}_{n}$, there is a simple closed
curve $\gamma$ in $\Omega$ that winds around each point in $\bigcup_{1\le
j\le n}\overline{W}_{j}$ and doesn't wind around any point outside $\Omega$.
(This is extremely easy to see if we first shrink each $\overline{W}_j$ to a
point.) Thus there is a univalent function $f$ on ${\mathbb{D}}$ such that $%
f({\mathbb{D}})$ is the set of points inside $\gamma$. As above we can
choose $s$, $0<s<1$, so that $\bigcup_{1\le j\le n}\overline{W}_{j}\subset
f(s{\mathbb{D}})$. It is clear that defining $\varphi(z)=f(z/s)$ and $r=1/s$
yields the desired function.
\end{proof}

\begin{lemma}
\label{c10} Suppose $S\in{\mathcal{B}}({\mathcal{H}})$, $r>1\ge r(S)$, and $%
\varphi$ is a univalent function on $r{\mathbb{D}}$. Then

\begin{enumerate}
\item $\sigma(\varphi(S))=\varphi(\sigma(S))$, $\sigma_{e}(\varphi(S))=%
\varphi(\sigma_{e}(S))$,

\item for every $\lambda\in\overline{{\mathbb{D}}}$, $\mathrm{{Ker}\,}%
(S-\lambda)=\mathrm{{Ker}\,}(\varphi(S)-\varphi(\lambda))$, $\mathrm{{Ker}\,}%
(S-\lambda)^{\ast}=\mathrm{{Ker}\,}(\varphi(S)-\varphi(\lambda))^{\ast}$, $%
\func{ind\,}(S-\lambda)=\func{ind\,}(\varphi(S)-\varphi(\lambda))$,

\item $\func{Lat\,}(S)=\func{Lat\,}(\varphi(S))$.
\end{enumerate}
\end{lemma}

\begin{proof}
Statement (1) is the spectral mapping theorem for the Riesz--Dunford
functional calculus. Since $\varphi$ is univalent, $\varphi(r{\mathbb{D}})$
is simply connected, and it follows from Runge's Theorem that $\varphi$ and $%
\varphi^{-1}$ are limits of polynomials that converge uniformly on compact
sets. This shows that (2) holds and then $S$ and $\varphi(S)$ generate the
same norm-closed unital algebras, from which (3) easily follows.
\end{proof}

Let ${\mathcal{Q}}({\mathcal{H}})$ denote the set of all operators in ${%
\mathcal{B}}({\mathcal{H}})$ whose nontrivial invariant subspaces all have
finite co-dimension.

\begin{lemma}
\label{c11}Suppose ${\mathcal{H}}$ is a separable infinite-dimensional
Hilbert space and $T\in{\mathcal{Q}}({\mathcal{H}})$. Then

\begin{enumerate}
\item $\sigma(T)$ is connected,

\item $\func{ind\,}(T-\lambda)\in\{0,-1\}$ for every $\lambda\in\rho_{S\text{%
-}F}(T)$.
\end{enumerate}
\end{lemma}

\begin{proof}
If $\sigma(T)$ were disconnected, then the Riesz--Dunford functional
calculus would give a spectral idempotent $P$ that commutes with $T$ and
either $\mathrm{{Ker}\,} P$ or $\mathrm{{Ker}\,}(1-P)$ would be a nonzero
invariant subspace with infinte co-dimension.

Clearly, $T$ cannot have any eigenvalues, and $T^{\ast}$ cannot have an
eigenvalue with an infinite-dimensional eigenspace. Hence, $\func{ind\,}%
(T-\lambda)\le 0$ whenever $\lambda\in\rho_{S\text{-}F}(T)$. Suppose that $%
\func{ind\,}(T-\lambda)\le -2$. Then $A=T-\lambda$ is injective, $A({%
\mathcal{H}})$ is closed and $\dim A({\mathcal{H}})^{\perp}\ge 2$. It
follows from the injectivity of $A$ that $\dim(A^{n}({\mathcal{H}})\ominus
A^{n+1}({\mathcal{H}}))\ge 2$ for every $n\ge 1$. Let $x_0,f\in A({\mathcal{H%
}})^{\perp}$ be two orthogonal unit vectors, and, for $n\ge 1$, take a unit
vector $x_{n}\in A^{n}({\mathcal{H}})\ominus A^{n+1}({\mathcal{H}})$ such
that $x_{n}\bot A^{n}f$. Let $M=\mathrm{{clos}\,}\mathrm{{span}\,}%
\{A^{n}f:n\ge 0\}$. Then $M$ is an invariant subspace for $T$ and if $%
x=\sum\alpha_{n}x_{n}\in M$, then, by successively projecting onto $A^{n}({%
\mathcal{H}})^{\perp}$, $n=1,2,\ldots$, we see that $x=0$. Hence $M$ has
infinite co-dimension, which is a contradiction to the assumption $T\in{%
\mathcal{Q}}$. Hence statement (2) holds.
\end{proof}

\begin{theorem}
\label{c12} The norm closure of the set ${\mathcal{R}}$ of operators with
connected essential spectrum and whose nontrivial invariant subspaces all
have finite co-dimension is the set ${\mathcal{T}}$ of all operators $T\in{%
\mathcal{B}}({\mathcal{H}})$, such that $\sigma(T)$ and $\sigma_{e}(T)$ are
both connected and such that 
\begin{equation*}
\func{ind\,}(T-\lambda)\in\{0,-1\} 
\end{equation*}
for every $\lambda\in\rho_{S\text{-}F}(T)$.
\end{theorem}

\begin{proof}
The inclusion $\mathrm{{clos}\,}{\mathcal{R}}\subset{\mathcal{T}}$ follows
from Lemma \ref{c11} and the fact that ${\mathcal{T}}$ is norm closed.
Suppose $T\in{\mathcal{T}}$ and $\varepsilon>0$. Let $N$ be a normal
operator such that $\sigma(N)=\sigma_{e}(N)=\sigma_{e}(T)$. It follows from
Proposition \ref{sim} that $T$ is in the closure of the similarity orbit of $%
T\oplus N$. We let $K=\sigma_{e}(T)$, and let $\mathcal{U}$ be the set of
those components of $\mathbb{C}\backslash K$ where $\func{ind\,}%
(T-\lambda)=-1$. Now we apply Lemma \ref{c9} to obtain $V,r,\varphi$, and
let $\gamma=\varphi({\mathbb{T}})$. Let $S$ be a quasianalytic shift
operator satisfying the conditions of Theorem~\ref{Shields} in Appendix. It
follows from Lemma \ref{c10} and Theorem \ref{Shields} that $\func{Lat\,} S=%
\func{Lat\,}\varphi(S)$, $\varphi(S)\in{\mathcal{R}}$, $\sigma(\varphi(S))=%
\varphi(\sigma(S))=\varphi(\overline{{\mathbb{D}}})$, $\sigma_{e}(%
\varphi(S))=\varphi(\sigma_{e}(S))=\varphi({\mathbb{T}})=\gamma$, $\func{%
ind\,}(\varphi(S)-\lambda)=-1$ exactly when $\lambda\in\varphi({\mathbb{D}})$%
.

Choose a normal operator $N_{\varepsilon}$ such that $\Vert
N-N_{\varepsilon}\Vert \le\varepsilon$ and $\sigma(N_{\varepsilon})=\mathrm{{%
clos}\,} [\sigma_{e}(T)_\varepsilon]$. It follows from Proposition \ref{sim}
that $T\oplus N_{\varepsilon}$ is in the closure of the similarity orbit of $%
\varphi(S)$; hence $T\oplus N_{\varepsilon}\in\mathrm{{clos}\,}{\mathcal{R}}$%
. Since $\Vert T\oplus N-T\oplus N_{\varepsilon}\Vert\le\varepsilon$ and $%
\varepsilon>0$ was arbitrary, we know that $T\oplus N\in\mathrm{{clos}\,}{%
\mathcal{R}}$. Since $T\in\mathrm{{clos}\,} {\mathcal{S}}(T\oplus N)$, we
know that $T\in\mathrm{{clos}\,}{\mathcal{R}}$. Thus, ${\mathcal{T}}\subset%
\mathrm{{clos}\,}{\mathcal{R}}$.
\end{proof}

\bigskip

We can now prove the main theorem of this section.

\begin{theorem}
\label{c13} Let $T\in{\mathcal{B}}({\mathcal{H}})$ and suppose that $%
\sigma_{e}(T)$ and $\sigma(T)$ are both connected. Then

\begin{enumerate}
\item If $\func{ind\,}(T-\lambda)\in\{0,-1\}$ for every $\lambda\in\rho_{S%
\text{-}F}(T)$ and if $0\ne M\in \func{Lat\,}_{ns}(T)$, then $\dim
M^{\perp}<\infty$ and $\sigma(T^{\ast}|M^{\bot})\subset\mathrm{{clos}\,}\{%
\overline{\lambda}\in\mathbb{C}: \func{ind\,}(T-\lambda)=-1\}$.

\item If $\func{ind\,}(T-\lambda)\in\{0,1\}$ for every $\lambda\in\rho_{S%
\text{-}F}(T)$ and if ${\mathcal{H}}\neq M\in\func{Lat\,}_{ns}(T)$, then $%
\dim M<\infty$ and $\sigma(T|M)\subset\mathrm{{clos}\,}\{\lambda\in\mathbb{C}%
: \func{ind\,}(T-\lambda)=1\}$.
\end{enumerate}
\end{theorem}

\begin{proof}
(1). In this case it follows from Theorem \ref{c12} that $T$ is a norm limit
of operators whose nonzero invariant subspaces all have finite co-dimension.
Since the set of projections with finite co-dimension is norm closed, it
follows that any nonzero norm-stable invariant subspace of $T$ has finite
co-dimension. If $U=\{\lambda\in\mathbb{C}:\func{ind\,}(T-\lambda)=-1\}
=\varnothing$, then it follows from Theorem \ref{c12} that $T$ and $T^{\ast}$
are norm limits of operators whose nontrivial invariant subspaces all have
finite co-dimension; whence, $\func{Lat\,}_{ns}(T)=\{0,1\}$. Suppose $%
U\neq\varnothing$. If $\varepsilon>0$, then $U$ has only finitely many
connected components $U_{1},\ldots,U_{m}$ containing a disc of radius $%
\varepsilon$. Since $\sigma(T)$ and $\sigma_{e}(T)$ are connected, the
components $U_{1},\ldots,U_{m}$ are simply connected. For each $k$, $1\le
k\le m$, there is an $r_k>1$ and a univalent mapping $\varphi_{k}:r_k{%
\mathbb{D}}\to U_k$ such that $U_{k}\backslash\varphi_{k}({\mathbb{D}}%
)\subset\{z\in U_{k}:\mathrm{{dist}\,}(z,\mathbb{C}\backslash
U_{k})<\varepsilon\}$. Let $S$ be a quasianalytic shift operator satisfying
the conditions of Theorem~\ref{Shields}, and let $A=\varphi_{1}(S)\oplus%
\cdots\oplus\varphi_{m}(S)\oplus N_{\varepsilon}$, where $N_{\varepsilon}$
is a normal operator with no eigenvalues whose spectrum is the closure of $%
\sigma_{e}(T)_{\varepsilon}$. It follows from Proposition \ref{sim} that $%
T\oplus N_{\varepsilon}$ is in the closure of the similarity orbit of $A$.
Also we can choose a normal operator $N$ with $\sigma(N)=\sigma_{e}(N)=%
\sigma_{e}(T)$ such that $\Vert N-N_{\varepsilon}\Vert <\varepsilon$. Since,
by Proposition \ref{sim}, $T$ is in the closure of the similarity orbit of $%
T\oplus N$, and $\sigma_{p}(A^{\ast})\subset\{\overline{\lambda}%
:\lambda\in\cup_{k=1}^{m}\varphi_{k}(\overline{{\mathbb{D}}})\}\subset \{%
\overline{\lambda}:\lambda\in U\}$, it follows that $T$ is a limit of
operators $T_{n}$ such that $\sigma_{p}(T_{n}^{\ast}|M_{n}^{\perp})\subset\{%
\overline{\lambda}:\lambda\in U\}$ for every $M_{n}\in\func{Lat\,} T_{n}$
with $\mathrm{{codim}\,}(M_{n})<\infty$. Thus statement (1) is proved.

(2). In this case $T^{\ast}$ satisfies the conditions of part (1), and since 
$Lat_{ns}(T^{\ast})$ is clearly $\{M^{\perp}:M\in\func{Lat\,}_{ns}(T)\}$,
the desired conclusion follows from part (1).
\end{proof}

\section{Lat$_{ns}$ for Weighted Unilateral Shifts}

We can use Theorem \ref{c13} to completely characterize $\func{Lat\,}_{ns}(T)
$ whenever $T$ is a weighted unilateral shift operator. The first step is
the following

\begin{corollary}
If $T$ is an injective unilateral weighted shift operator, then every
nonzero norm-stable invariant subspace of $T$ has finite co-dimension.
\end{corollary}

\begin{proof}
If the weights of $T$ are not bounded away from $0$, then a compact
perturbation of $T$, obtained by replacing a subsequence of weights with
zeros, is quasidiagonal, and hence biquasitriangular. Since weighted shifts
have connected spectrum, it follows from \cite[Corollary 3.9]{{AFS}} that $T$
has no nontrivial norm-stable invariant subspaces. If the weights of $T$ are
bounded away from $0$, then $T$ satisfies condition (1) in Theorem \ref{c13}.
\end{proof}

The following lemma is a simple application of the Gram-Schmidt process.

\begin{lemma}
\label{basis}Suppose $\{u_1,\ldots,u_n\}$ is a linear basis for a subspace $M
$ of ${\mathcal{H}}$, and suppose, for $k\in\mathbb{N}$, that $M_k$ is the
span of $\{u_{k1},\ldots,u_{kn}\}$. If $\lim_{k\to\infty}\Vert
u_{j}-u_{kj}\Vert=0$ for $1\le j\le n$, then $\Vert P_{M_{k}}-P_{M}\Vert\to 0
$.
\end{lemma}

\begin{corollary}
If $T\in{\mathcal{B}}({\mathcal{H}})$ and $M_1,\ldots,M_k$ are
one-dimensional subspaces in $\func{Lat\,}_{ns}(T)$ \textrm{(}respectively, $%
\func{Lat\,}_{s}(T)$\textrm{)}, then $M_{1}+\ldots+M_{k}$ belongs to $\func{%
Lat\,}_{ns}(T)$ \textrm{(}respectively, $\func{Lat\,}_{s}(T)$\textrm{)}.
\end{corollary}

\begin{lemma}
\label{c15} Suppose $T$ is a bounded operator on ${\mathcal{H}}$, and $M$ is
a nonzero finite-dimensional cyclic invariant subspace for $T^{\ast}$ such
that $\sigma(T^{\ast}|M)\cap\sigma_{e}(T^{\ast})=\varnothing$ and $%
\sigma(T^{\ast}|M)\cap\{\overline{\lambda}:\lambda\in\sigma_{p}(T)\}=%
\varnothing$. Then $M\in\func{Lat\,}_{ns}(T^{\ast})$.
\end{lemma}

\begin{proof}
We know that if $p(z)$ is the minimal polynomial for $T^{\ast}|M$, then its
set of roots is $\sigma(T^{\ast}|M)$ and $M\subset\mathrm{{Ker}\,}%
(p(T^{\ast}))$. Since $\sigma(T^{\ast}|M)\cap\{\overline{\lambda}%
:\lambda\in\sigma_{p}(T)\}=\varnothing$, we know that $\mathrm{{Ker}\,}%
(p(T^{\ast})^{\ast})=0$. Moreover, $\sigma(T^{\ast}|M)\cap\sigma_{e}(T^{%
\ast})=\varnothing$ implies that $p(T^{\ast})$ is Fredholm, and we know that 
$\mathrm{{Range}\,}(p(T^{\ast}))=\mathrm{{Ker}\,}(p(T^{\ast})^{\ast})
^{\perp}={\mathcal{H}}$, thus $p(T^{\ast})$ is surjective. If $\{A_{k}\}$ is
a sequence converging in norm to $T^{\ast}$, then $\Vert
p(A_{k})-p(T^{\ast})\Vert\to 0$, and, by \cite[Lemma 1.6]{CH}, the
projection $Q_{k}$ onto $\mathrm{{Ker}\,}(p(A_k))$ converge in norm to the
projection $Q$ onto $\mathrm{{Ker}\,}(p(T^{\ast}))$. Let $e$ be a cyclic
vector for $T^{\ast}|M$. If $m=\deg(p)$, then $\{e,T^{\ast}e,\ldots,(T^{%
\ast})^{m-1}e\}$ is a basis for $M$, and since $Q_{k}e\in\mathrm{{Ker}\,}
p(A_k)$, the set $\{Q_{k}e,A_{k}Q_{k}e,\ldots,A_{k}^{m-1}Q_{k}e\}$ spans a
subspace $M_{k}\in\func{Lat\,}(A_{k})$ for each $k\ge 1$. It follows from
Lemma \ref{basis} that $\Vert P_{M_{k}}-P_{M}\Vert\to 0$. Thus $M\in\func{%
Lat\,}_{ns}(T^{\ast})$.
\end{proof}

We are now ready to completely characterize $\func{Lat\,}_{ns}(T)$ for every
unilateral shift operator $T$.

\begin{theorem}
\label{wtshift} Suppose $T$ is a weighted unilateral shift operator and let $%
r_0=\inf\{\vert\lambda\vert:\lambda\in \sigma_{e}(T)\}$. If $r_0=0$, then $%
\func{Lat\,}_{ns}(T)=\{0,1\}$. If $r_0>0$, then a subspace $M\neq{\mathcal{H}%
}$ is in $\func{Lat\,}_{ns}(T^{\ast})$ if and only if $\dim M<\infty$ and $%
\vert\lambda\vert \le r_{0}$ for every $\lambda\in\sigma(T^{\ast}|M)$.
\end{theorem}

\begin{proof}
The "only if" part follows from Theorem \ref{c13}. Suppose that \newline
\noindent $\{e_0,e_1,\ldots\}$ is an orthonormal basis, and $%
\{\alpha_n\}_{n\ge 0}$ is a sequence of positive numbers such that $%
Te_n=\alpha_n e_{n+1}$ for all $n\ge 0$. Then $T^{\ast}e_0=0$ and $%
T^{\ast}e_{n+1}=\alpha_n e_n$ for all $n\ge 0$. Suppose that $\lambda\in%
\mathbb{C}$ and $\lambda\in\sigma_{p}(T^{\ast})$. Then there is a vector $%
0\neq f_{\lambda,1}=(\beta_0,\beta_1,\ldots)$ such that $(T^{\ast}-%
\lambda)f_{\lambda,1}=0$. It is clear that $\beta_0\neq 0$, so we can assume 
$\beta_0=1$, and, for $n\ge 1$, 
\begin{equation*}
\beta_n=\frac{\lambda^n}{\alpha_{0}\cdot\ldots\cdot\alpha_{n-1}}. 
\end{equation*}
It follows that $\mathrm{{Ker}\,}(T^{\ast}-\lambda)$ is $1$-dimensional.
Thus, if $M$ is a finite-dimen\-sional invariant subspace for $T^{\ast}$,
then each eigenvalue for $T^{\ast}|M$ has exactly one Jordan block in its
Jordan form. Thus $T^{\ast}|M$ is cyclic. It follows from Lemma \ref{c15}
that if $\sigma(T^{\ast}|M)\cap\sigma_{e}(T^{\ast})=\varnothing$, then $M\in%
\func{Lat\,}_{ns}(T^{\ast})$. Next suppose that there exists an $%
f_{\lambda,2}$ such that $(T^{\ast}-\lambda)f_{\lambda,2}=f_{\lambda,1}$.
Then we can choose $f_{\lambda,2}=(0,\gamma_1,\gamma_2,\ldots)$, and we see
the $\gamma_{k}$'s are uniquely determined and, for $n\ge 1$, 
\begin{equation*}
\gamma_{n}=\frac{n\lambda^{n-1}}{\alpha_{0}\cdot\ldots\cdot\alpha_{n-1}}. 
\end{equation*}
More generally, if we have $f_{\lambda,1},\ldots,f_{\lambda,m}$ such that,
for $1\le k<m$, 
\begin{equation*}
(T^{\ast}-\lambda)f_{\lambda,k+1}=f_{\lambda,k} 
\end{equation*}
and the first $k$ coordinates of $f_{\lambda,k+1}$ are $0$, then 
\begin{equation*}
f_{\lambda,k}=(0,\ldots,0,c_{k,k},c_{k,k+1}\lambda,c_{k,k+2}\lambda^{2},%
\ldots) , 
\end{equation*}
where the positive numbers $c_{k,j}$ ($j\ge k$) depend only on the weights $%
\{\alpha_n\}$ and not on $\lambda$. Note that $\Vert f_{\lambda,k}\Vert$
depends on $\vert\lambda\vert$, so that if $\Vert f_{\lambda_{0},k}\Vert
<\infty$ for some $\lambda_0$, then $\Vert f_{\lambda,k}\Vert<\infty$ for
all $\lambda$ with $\vert\lambda\vert\le\vert\lambda_0\vert=r$. Moreover,
the map $\lambda\mapsto\Vert f_{\lambda,k}\Vert$ is continuous on $\overline{%
r{\mathbb{D}}}$ (by the dominated convergence theorem). If a sequence $%
\{h_n\}$ of vectors in a Hilbert space converges weakly to $h$ and if $\Vert
h_n\Vert\to\Vert h\Vert$, then $\Vert h_{n}-h\Vert \to 0$. It follows that
the map $\lambda\mapsto f_{\lambda,k}$ is norm continuous on $\overline{r{%
\mathbb{D}}}$.

Now suppose $M\in\func{Lat\,}(T^{\ast})$ is finite-dimensional and $\vert
\lambda\vert \le r_0$ for every $\lambda\in\sigma(T^{\ast}|M)$. If $%
p(z)=(z-\lambda_1)^{m_{1}}\cdot\ldots\cdot(z-\lambda_s)^{m_{s}}$ is the
minimal polynomial for $T^{\ast}|M$, then $\{f_{\lambda_{j},k}:1\le j\le
s,1\le k\le m_{j}\}$ is a linear basis for $M$. The desired conclusion
follows from Lemma \ref{basis}.
\end{proof}

We see that quasianalytic shifts are points of norm continuity of $\func{%
Lat\,}$.

\begin{corollary}
Suppose $T$ is a quasianalytic unilateral shift operator satisfying the
conditions of Theorem~\ref{Shields}. Then $\func{Lat\,}(T)=\func{Lat\,}%
_{ns}(T)$.
\end{corollary}

\section{A Result for Lat$_{s}$.}

Suppose that $T\in {\mathcal{B}}({\mathcal{H}})$ and $P\in\func{Lat\,}(T)$.
We define the \emph{index of} $T$ \emph{relative to} $P$ by 
\begin{equation*}
\func{ind\,}(T,P)=\dim P({\mathcal{H}})\ominus \mathrm{{clos}\,} TP({%
\mathcal{H}}). 
\end{equation*}
If $T$ is bounded from below, i.e., $\mathrm{{Ker}\,} T=0$ and $T({\mathcal{H%
}})$ is closed, then $\func{ind\,}(T,P)$ can be defined in terms of the
Fredholm index, namely $\func{ind\,}(T,P)=-\func{ind\,}(T|P)$. We first
prove a semicontinuity result.

\begin{lemma}
\label{semicontinuous} Suppose that $T,S_1,S_2,\ldots\in {\mathcal{B}}({%
\mathcal{H}})$ and $P\in\func{Lat\,} T$ and $P_{n}\in \func{Lat\,} S_n$ for $%
n\ge 1$ are such that

\begin{enumerate}
\item $T$ is bounded from below,

\item $\Vert S_n-T\Vert \to 0$, and

\item $P_n\to P$ in the strong operator topology.
\end{enumerate}

Then 
\begin{equation*}
\func{ind\,}(T,P)\le\liminf_{n\to\infty}\func{ind\,}(S_n,P_n). 
\end{equation*}
\end{lemma}

\begin{proof}
Let $\func{ind\,}(T,P)\ge k$, and let $E$ be a linear subspace of $P({%
\mathcal{H}})\ominus TP({\mathcal{H}})$, $\dim E=k$. Since $\dim E<\infty$, $%
P_n|E$ converge to $1|E$ in norm.

Suppose that $P_nE\cap S_nP_n({\mathcal{H}})\neq\{0\}$ for large $n$.
Passing to a subsequence, we can find $e_n\in E$, $\|e_n\|=1$, such that 
\begin{equation*}
P_ne_n\in S_nP_n ({\mathcal{H}}),\qquad e_n\to e\in E, 
\end{equation*}
and hence $Pe=e\neq 0$. Denote $u_n=S_n^{-1}P_ne_n\in P_n ({\mathcal{H}})$.
Then $\{u_n\}$ is a bounded sequence, and then, a Cauchy sequence (since $%
S_nu_n\to e$, $n\to\infty$); denote $f=\lim_{n\to\infty}u_n$. We have $%
Tf=e\in E$, $f\not\in P({\mathcal{H}})$. Furthermore, $u_n=P_nu_n\to Pf\neq f
$. This contradiction shows that $P_nE\cap S_nP_n ({\mathcal{H}})=\{0\}$ for
large $n$, and, hence, 
\begin{equation*}
\dim P_n({\mathcal{H}})\ominus S_nP_n({\mathcal{H}})\ge \dim E=k. 
\end{equation*}
\end{proof}

\begin{lemma}
\label{unweighted2} Suppose that $S$ is the (unweighted) unilateral shift
operator, $r>1$, and $\varphi:r{\mathbb{D}}\rightarrow\mathbb{C}$ is
univalent with $\varphi(0)=0$. Then for every $P\in \func{Lat\,} S$ we have $%
\func{ind\,}(\varphi(S),P)\le 1$.
\end{lemma}

\begin{proof}
It follows from Lemma \ref{c10} that $\func{Lat\,} S=\func{Lat\,}\varphi(S)$%
. Since $\varphi$ is univalent and $\varphi(0)=0$, we have $%
\varphi(z)=z\psi(z)$ where $\psi(z)\neq 0$ for every $z\in r{\mathbb{D}}$.
Hence $\varphi(S)=SA$ with $A$ invertible and $A$ and $A^{-1}=(\frac{1}{\psi}%
)(S)$ are in the weakly closed algebra generated by $S$. Thus $\varphi(S)P({%
\mathcal{H}})=SP({\mathcal{H}})$, and we conclude that $\func{ind\,}%
(\varphi(S),P)=\func{ind\,}(S,P)\le 1$. The last inequality follows from
Beurling's characterization of $\func{Lat\,} S$ \cite{Beu}.
\end{proof}

\begin{theorem}
\label{stronglystable} Suppose that $T\in {\mathcal{B}}({\mathcal{H}})$, $%
\sigma(T)$ and $\sigma_{e}(T)$ are both connected, and $\func{ind\,}%
(T-\lambda)\in\{-1,0\}$ for every $\lambda\in\rho_{S\text{-}F}(T)$, and
suppose that $\func{ind\,}(T)=-1$ and $\mathrm{{Ker}\,}(T)=0$. If $P\in 
\func{Lat\,}_{s}(T)$, then $\func{ind\,}(T,P)\le 1$.
\end{theorem}

\begin{proof}
We imitate the proof of Theorem \ref{c12} replacing a quasianalytic shift
with the unweighted unilateral shift $S$. Since $\func{ind\,}(T)=-1$, we can
assume (by composing with a disc automorphism) that $\varphi(0)=0$. It
follows that $T$ is the norm limit of a sequence $\{S_n\}$, where each $S_n$
is similar to $\varphi_n(S)$ for some univalent functions $\varphi_n$ on a
neighborhood of $\overline{{\mathbb{D}}}$ with $\varphi_{n}(0)=0$. The
desired conclusion follows from Lemmas \ref{semicontinuous} and \ref%
{unweighted2}.
\end{proof}

Suppose $M\in \func{Lat\,}(T)$ and $N=\mathrm{{clos}\,} T(M)$. With respect
to the decomposition ${\mathcal{H}}=N\oplus(M\ominus N)\oplus M^{\bot}$, $T$
has an operator matrix 
\begin{equation*}
\left( 
\begin{array}{ccc}
A & B & C \\ 
0 & 0 & D \\ 
0 & 0 & E%
\end{array}
\right) 
\end{equation*}
where the size of the $0$ in the $(2,2)$-entry is $\func{ind\,}(T,M)$. In 
\cite{ABFP} C.~Apostol, H.~Bercovici, C.~Foia\c{s}, and C.~Pearcy introduced
a class $\mathbb{A}_{\aleph_{0}}$ of contraction operators and they proved
that if $T\in\mathbb{A}_{\aleph_{0}}$, then $T$ has invariant subspaces with
arbitrary index. They also proved that if $T$ is a contractive unilateral
weighted shift whose weights converge to $1$, then either $T\in\mathbb{A}%
_{\aleph_{0}}$ or $T$ is similar to the unweighted unilateral shift. Thus
the Bergman shift with weights $\sqrt{\frac{n+1}{n+2}}$ is in $\mathbb{A}%
_{\aleph_{0}}$.

\begin{corollary}
Suppose $T$ is a unilateral shift whose weights converge to $\Vert T\Vert
\neq 0$. Then the following are equivalent:

\begin{enumerate}
\item $\func{Lat\,}(T)=\func{Lat\,}_{s}(T)$

\item $T$ is similar to $\Vert T\Vert $ times the unweighted unilateral
shift,

\item $\func{ind\,}(T,P)\le 1$ for every $P\in \func{Lat\,}(T)$.
\end{enumerate}
\end{corollary}

\begin{remark}
If, in the proof of Theorem \ref{stronglystable}, we replace the role of the
unweighted unilateral shift $S$ with a direct sum $S^{(n)}$ of $n\ge 1$
copies of $S$, then we can show that if $\func{ind\,}(T)=-n$, $\mathrm{{Ker}%
\,}(T)=0$, $\func{ind\,}(T-\lambda)\in\{-n,0\}$ for every $\lambda\in\rho_{S%
\text{-}F}(T)$, and $\sigma(T)$ and $\sigma_{e}(T)$ are both connected, then 
$\func{ind\,}(T,P)\le n$ for every $P\in \func{Lat\,}_{s}(T)$.
\end{remark}

\section{Questions}

We conclude with some open questions.

\textbf{Question }$1$. Suppose that $T\in{\mathcal{Q}}({\mathcal{H}})$,
i.e., the nontrivial invariant subspaces of $T$ all have finite
co-dimension, must $\sigma_{e}(T)$ be connected? In particular, is there an
operator $T\in{\mathcal{Q}}({\mathcal{H}})$ whose spectrum is an annulus and
essential spectrum is its boundary with Fredholm index $-1$ inside the
annulus?

\textbf{Question }$2.$ Suppose that $T$ is a weighted unilateral shift with
closed range. If $P\in\func{Lat\,}(T)$ and $\func{ind\,}(T,P)=1$, must $P\in 
\func{Lat\,}_{s}(T)$?

\textbf{Question }$3$. What is $\func{Lat\,}_{s}(T)$ when $T$ is a
unilateral weighted shift? What if $T$ is the Bergman shift?

\appendix

A key ingredient of the proof of our main results on stable invariant
subspaces involves properties of \emph{quasianalytic shift operators}, which
are weighted unilateral shifts with weights converging to $1$, whose
essential spectrum is the unit circle, whose spectrum is the closed unit
disc, whose Fredholm index is $-1$ on the open unit disc. An example of a
quasianalytic shift has weights $\exp(\sqrt{n+1}-\sqrt{n})$. These shifts
were used in \cite{HNRR} and \cite{H} to show that results of Lomonosov \cite%
{L1}, \cite{L2} did not lead to an immediate solution of the invariant
subspace problem. The most important property of quasianalytic shifts from
our point of view concerns their invariant subspaces. Question 17 in the
seminal 1974 paper \cite{Shi} on weighted shift operators by Allen Shields
is whether the nonzero invariant subspaces of quasianalytic shifts all have
finite co-dimension. In this section we establish a version of Domar's
result answering Shields' question in the affirmative.

Let us introduce some definitions. Given a function $\omega:\mathbb{Z}%
_{+}\rightarrow\lbrack1,+\infty)$ such that 
\begin{gather*}
0<\inf_{n\ge 0}\frac{\omega(n+1)}{\omega(n)}\le\sup_{n\ge 0}\frac {%
\omega(n+1)}{\omega(n)}<\infty, \\
\lim_{n\rightarrow\infty}\omega(n)^{1/n}=1,
\end{gather*}
we consider the unilateral shift operator $S_{\omega}:\ell^{2}(\mathbb{Z}%
_{+})\rightarrow\ell^{2}(\mathbb{Z}_{+})$, $S_{\omega}e_{n}=(\omega
(n+1)/\omega(n))e_{n+1}$, where $e_{n}=\{\delta_{mn}\}_{m\ge 0}\in\ell ^{2}(%
\mathbb{Z}_{+})$.

Let $A_{\omega}^{2}$ be the Beurling space of the functions $f$ analytic in
the unit disc $\mathbb{D}$, with $f(z)=\sum_{n\ge 0}\hat{f}(n)z^{n}$, $z\in%
\mathbb{D}$, such that 
\begin{equation*}
\Vert f\Vert_{A_{\omega}^{2}}^{2}=\sum_{n\ge 0}|\hat{f}(n)|^{2}\omega
(n)^{2}<\infty. 
\end{equation*}
The operator $M_{z}:f\mapsto zf$ of multiplication by the independent
variable on $A_{\omega}^{2}$ is isomorphic to $S_{\omega}$. Denote $\omega
_{s}(n)=\omega(n)(1+n)^{-s}$, and suppose that the sequence $\log\omega
_{1}(n)$ is convex for large $n$. Then $\omega_{1}(n+m)\le c\,\omega
_{1}(n)\omega_{1}(m)$, $n\ge 0$, $m\ge 0$, and $A_{\omega}^{2}$ is a Banach
algebra with respect to the convolution multiplication: 
\begin{gather*}
\Vert fg\Vert_{A_{\omega}^{2}}^{2}=\sum_{n\ge 0}|\widehat{fg}%
(n)|^{2}\omega(n)^{2}=\sum_{n\ge 0}\bigl|\sum_{0\le k\le n}\hat{f}(k)\hat {g}%
(n-k)\bigr|^{2}\omega(n)^{2}\le \\
\sum_{n\ge 0}\Bigl[\sum_{0\le k\le n}\Bigl(\frac{\omega(n)}{\omega
(k)\omega(n-k)}\Bigr)^{2}\Bigr]\times \\
\Bigl[\sum_{0\le k\le n}|\hat{f}(k)|^{2}\omega(k)^{2}|\hat{g}%
(n-k)|^{2}\omega(n-k)^{2}\Bigr]\le \\
c\sum_{m\ge 0}|\hat{f}(m)|^{2}\omega(m)^{2}\cdot\sum_{v\ge 0}|\hat{g}%
(v)|^{2}\omega(v)^{2}\times \\
\max_{n\ge 0}\sum_{0\le k\le n}\frac{(n+1)^{2}}{(k+1)^{2}(n-k+1)^{2}}\le \\
c\Vert f\Vert_{A_{\omega}^{2}}^{2}\Vert g\Vert_{A_{\omega}^{2}}^{2}.
\end{gather*}
By \cite[Corollary 1, p.94]{Shi}, the space of maximal ideals of $A_{\omega
}^{2}$ is the closed unit disc $\bar{\mathbb{D}}$.

Next, we suppose that $\omega$ is log-convex, that is it extends to $\mathbb{%
R}_{+}$ in such a way that $t\mapsto\log\omega(\exp t)$ is convex for large $%
t$, and that 
\begin{equation*}
\sum_{n\ge 0}\frac{\log\omega(n)}{n^{3/2}+1}=\infty. 
\end{equation*}
In this case, the space $A_{\omega}^{2}$ is quasi-analytic, in the sense
that $A_{\omega}^{2}\subset C^{\infty}(\bar{\mathbb{D}})$, and if $f\in
A_{\omega }^{2}$, $z\in\bar{\mathbb{D}}$, $f^{(n)}(z)=0$ for all $n\ge 0$,
then $f=0$, see \cite{Ca,Sa,Ko}.

Let us assume that $\omega$ satisfies an additional property: for every
polynomial $p$, there exists $c(p)<\infty$ such that 
\begin{equation}
\Vert pf_{1}f_{2}\Vert_{A_{\omega}^{2}}\le c(p)\Vert pf_{1}\Vert_{A_{\omega
}^{2}}\Vert pf_{2}\Vert_{A_{\omega}^{2}},\qquad f_{1},f_{2}\in
A_{\omega}^{2}.   \label{wa}
\end{equation}
Then by Theorem 1 of Domar \cite{D}, every non-trivial closed $M_{z}$%
-invariant subspace of $A_{\omega}^{2}$ has finite co-dimension.

By Remarks 1 and 2 of \cite{D}, we need only to verify (\ref{wa}) for $%
p(z)=(z-\lambda)^{n}$ with $\lambda\in\bar{\mathbb{D}}$. In fact, Domar
gives (Theorem 3, Remarks 3 and 4 in \cite{D}, see also \cite{DMZ}) a
sufficient condition on $\omega$ in order that the space 
\begin{equation*}
A_{\omega}^{1}=\bigl\{\sum_{n\ge 0}\hat{f}(n)z^{n}:\sum_{n\ge 0}|\hat {f}%
(n)|\omega(n)<\infty\bigr\}
\end{equation*}
would satisfy the property that every non-trivial closed $M_{z}$-invariant
subspace of $A_{\omega}^{1}$ has finite co-dimension.

Here we present two arguments giving (\ref{wa}). For simplicity, let us
assume that $p(z)=z-1$. We need to prove that 
\begin{equation}
\|(z-1)f_{1}f_{2}\|_{A^{2}_{\omega}}\le
c\|(z-1)f_{1}\|_{A^{2}_{\omega}}\cdot\|(z-1)f_{2}\|_{A^{2}_{\omega}},\qquad
f_{1},f_{2}\in A^{2}_{\omega }.   \label{wb}
\end{equation}

The first argument. It is known (see, for example, \cite[Lemma 5.2]{Es}
adapting a general construction in \cite[Proposition B.1]{BH}) that there
exists a positive continuous function $\varphi:(1,2]\rightarrow(0,\infty)$
such that 
\begin{equation*}
\int_{1}^{2}\varphi\left( r\right) r^{-2n}dr\asymp\frac{1}{\omega(n)},\qquad
n\ge 0. 
\end{equation*}
Let $\Omega=\{z:1<|z|<2\}$. Given a function $g$ in 
\begin{equation*}
L^{2}(\varphi,\Omega)=\{h:\int_{\Omega}|h|^{2}\varphi^{2}<\infty\}, 
\end{equation*}
we define 
\begin{equation*}
\mathcal{C}g(z)=\frac{1}{\pi}\int_{1<|w|<2}\frac{g(w)}{z-w}\,dm_{2}(w). 
\end{equation*}
Then the space $A_{\omega}^{2}$ coincides (\cite{Dy,BH}) with the image
space $\{\mathcal{C}g:g\in L^{2}(\varphi,\Omega)\}$, 
\begin{equation*}
\Vert f\Vert_{A_{\omega}^{2}}\asymp\inf\{\Vert g\Vert_{L^{2}(\varphi,\Omega
)}:\mathcal{C}g|\mathbb{D}=f\}. 
\end{equation*}
Now, let $g_{1}(z)=(z-1)f_{1}(z)$, $g_{2}(z)=(z-1)f_{2}(z)$, where $%
g_{1},g_{2}\in L^{2}(\varphi,\Omega)$ are canonical densities (\cite[Section
2]{BH}). Then 
\begin{equation*}
\frac{\mathcal{C}g_{1}\cdot\mathcal{C}g_{2}}{\cdot-1}=\mathcal{C}\left( \bar{%
\partial}\left( \frac{\mathcal{C}g_{1}\cdot\mathcal{C}g_{2}}{\cdot -1}%
\right) \right) , 
\end{equation*}
and to prove (\ref{wb}) we need only to verify that 
\begin{equation*}
\bigl\|g_{1}\cdot\frac{\mathcal{C}g_{2}}{\cdot-1}\bigr\|_{L^{2}(\varphi
,\Omega)}\le c\Vert g_{1}\Vert_{L^{2}(\varphi,\Omega)}\cdot\Vert
g_{2}\Vert_{A_{\omega}^{2}}. 
\end{equation*}
Finally, the equality 
\begin{equation*}
\frac{\left( \mathcal{C}g_{2}\right) (z)}{z-1}=\frac{\mathcal{C}g_{2}(z)-%
\mathcal{C}g_{2}(1)}{z-1}=\left( \mathcal{C}(\frac{g_{2}}{\cdot -1})\right)
(z), 
\end{equation*}
Lemma 8.2, Proposition 8.4, and Lemma B.6 in \cite{BH} give us the estimate 
\begin{equation*}
\bigl\|\frac{\mathcal{C}g_{2}}{\cdot-1}\bigr\|_{L^{\infty}(\Omega)}\le c%
\bigl\|\frac{g_{2}}{\cdot-1}\bigr\|_{L^{3}(\Omega)}, 
\end{equation*}
if $\log\omega_{3}(n)$ is concave for large $n$. \bigskip

The second argument was proposed to us by Omar El-Fallah.

First of all, if $\omega_{2}(n)$ increases for large $n$, then 
\begin{equation}
\Vert(z-1)f\Vert_{A_{\omega}^{2}}\ge c\Vert
f\Vert_{A_{\omega_{1}}^{2}},\qquad f\in A_{\omega}^{2}.   \label{wc}
\end{equation}
Indeed, if $g(z)=(z-1)f(z)$, $g(z)=\sum_{n\ge 0}\hat{g}(n)z^{n}$, then 
\begin{gather*}
f(z)=\frac{g(z)-g(1)}{z-1}=\sum_{0\le k<n}\hat{g}(n)z^{k}, \\
\Vert f\Vert_{A_{\omega_{1}}^{2}}^{2}=\sum_{n\ge 0}|\hat{f}(n)|^{2}\omega
_{1}(n)^{2}=\sum_{n\ge 0}\bigl|\sum_{k>n}\hat{g}(k)\bigr|^{2}%
\omega_{1}(n)^{2}\le \\
c\sum_{n\ge 0}\frac{\omega_{1}(n)^{2}}{n+1}\sum_{k>n}|\hat{g}%
(k)|^{2}k^{2}=c\sum_{k>0}|\hat{g}(k)|^{2}k^{2}\sum_{0\le n<k}\frac{%
\omega_{1}(n)^{2}}{n+1}\le \\
c\sum_{k>0}|\hat{g}(k)|^{2}\omega(k)^{2}=c\Vert g\Vert_{A_{\omega}^{2}}^{2}.
\end{gather*}
Next, we note that 
\begin{equation*}
\Vert f\Vert_{A_{\omega}^{2}}\asymp|f(0)|+\Vert f^{\prime}\Vert_{A_{\omega
_{1}}^{2}}^{2},\qquad f\in A_{\omega}^{2}. 
\end{equation*}
Therefore, by (\ref{wc}), to prove (\ref{wb}) we need only to verify that 
\begin{equation*}
\Vert\lbrack(z-1)f]^{\prime}g\Vert_{A_{\omega_{1}}^{2}}^{2}\le c\Vert
\lbrack(z-1)f]^{\prime}\Vert_{A_{\omega_{1}}^{2}}^{2}\cdot\Vert g\Vert
_{A_{\omega_{1}}^{2}}^{2};
\end{equation*}
this estimate follows immediately if $\log\omega_{2}(n)$ is concave for
large $n$. \bigskip

Summing up, we obtain the following result.

\begin{theorem}
\label{Shields} Let $\omega:\mathbb{Z}_{+}\rightarrow\lbrack1,+\infty)$ be
such that $\lim_{n\rightarrow\infty}\omega(n)^{1/n}=1$, for every $k\ge 0$,
the sequence $\omega(n)(n+1)^{-k}$ is concave for large $n$, $\omega$ is
log-convex, and 
\begin{equation*}
\sum_{n\ge 0}\frac{\log\omega(n)}{n^{3/2}+1}=\infty. 
\end{equation*}
Then every nontrivial closed invariant subspace of $S_{\omega}$ has finite
co-dimen\-sion.
\end{theorem}

\end{document}